\documentclass[a4paper,DIV=12,10pt]{scrartcl}

\usepackage[utf8]{inputenc}
\usepackage{amsfonts,amssymb,amsmath}
\usepackage{mathtools}

\usepackage{epsfig}
\usepackage{MnSymbol}
\usepackage{booktabs}
\usepackage{cite}
\usepackage{float}
\usepackage{hyperref}
\usepackage{siunitx}

\usepackage{wrapfig}

\usepackage{graphics}
\usepackage{color}

\numberwithin{equation}{section}
\numberwithin{table}{section}
\numberwithin{figure}{section}

\newcommand{\ud}{\,\mathrm{d}}
\newcommand{\R}{\mathbb{R}}
\newcommand{\N}{\mathbb{N}}

\newcommand{\mvec}{\boldsymbol}

\numberwithin{equation}{section}
\numberwithin{figure}{section}
\numberwithin{table}{section}

\newtheorem{defi}{Definition}[section]
\newtheorem{problem}[defi]{Problem}
\newtheorem{rem}[defi]{Remark}
\newtheorem{theorem}[defi]{Theorem}
\newtheorem{lemma}[defi]{Lemma}
\newenvironment{mproof}{\paragraph{Proof}}{\hfill$\blacksquare$}

\usepackage{cleveref}

\usepackage{caption}
\usepackage{subcaption}

\begin{document}
	
	\title{Optimal order FEM for dynamic poroelasticity: Error analysis for equal order elements }
	
	\author{
		Markus Bause$^\dag$\thanks{bause@hsu-hh.de (corresponding author)}\;, Mathias Anselmann$^\dag$\\
		{\small $^\dag$ Helmut Schmidt University, Faculty of
			Mechanical and Civil Engineering, Holstenhofweg 85,}\\ 
		{\small 22043 Hamburg, Germany}
	}
	\date{}
	
	\maketitle

\begin{abstract}
	The numerical approximation of dynamic poroelasticity, modeling flow in deformable porous media, by a family of continuous space-time finite element methods is investigated. Equal order approximation in space without any further stabilization is used for the displacement and pore pressure variable. Optimal order $L^\infty(L^2)$  error estimates are proved and numerically confirmed.
\end{abstract}

\section{Introduction}
\label{Sec:Intro}

We study the numerical approximation by equal order finite element methods of dynamic poroelasticity in the open bounded domain $\Omega \subset \R^d$, with $d\in \{2,3\}$, and the time interval $I=(0,T]$, with $T>0$, written as a first-order sytem in time
\begin{subequations}
\label{Eq:HPS}
\begin{alignat}{4}
 \partial_t \mvec u - \mvec v & = 0\,,  &   
\rho \partial_t \mvec v - \nabla \cdot (\mvec C \mvec \varepsilon (\mvec u)) + \alpha \mvec \nabla p  & = \rho \mvec f\\
&& c_0\partial_t p + \alpha \nabla \cdot \partial_t \mvec u  -  \nabla \cdot (\mvec K \nabla p)  & = g
\end{alignat}
\end{subequations}
and equipped with the initial conditions $	\mvec u (0) = \mvec u_0$, $\mvec v (0) = \mvec v_0$ and $p(0) = p_0$ and homogeneous Dirichlet boundary conditions for $\mvec u$ and $p$.  The quantity $\mvec \varepsilon (\mvec u):= (\nabla \mvec u + (\nabla \mvec u)^\top)/2$ is the symmetrized gradient. The coefficients $\rho>0$, $\alpha>0$ and $c_0 >0$ as well as the symmetric and positive definite tensors $\mvec C$ and $\mvec K$ are  constants. Well-posedness of \eqref{Eq:HPS} follows from \cite{MB_JR18,MB_STW22,MB_S89}. Its energy is measured by \mbox{\cite[p.~15]{MB_J16}}
\begin{equation}
	\label{Eq:HPSenergy}
	E(t) := \int_I \left(\frac{\rho}{2} \langle \partial_t \mvec u(t), \partial_t \mvec u(t)\rangle + \langle \mvec C \mvec \varepsilon (t), \mvec \varepsilon(t)\rangle + \frac{c_0}{2} \langle p(t),p(t)\rangle \right) \ud t\,. 
\end{equation}

 In \cite{MB_BAKR22}, a continuous Galerkin method in space and time, that is based on an approximation by piecewise polynomials of order $k$ in time and the Taylor--Hood family of finite element spaces (cf., e.g., \cite{MB_BF91}) with piecewise polynomials of order $r+1$ for $\mvec u$ and $\mvec v$ as well as piecewise polynomials of order $r$ for $p$ in space, is investigated thoroughly. By a generalization of the energy techniques used in \cite{MB_KM04}, for the approximation $(\mvec u_{\tau},\mvec v_{\tau,h}, p_{\tau,h})$ of $(\mvec u,\mvec v,p)$ the error estimate 
\begin{equation*}
	\label{MainEstEng}
	\begin{aligned}
	 \|\nabla ( \mvec u(t) -\mvec u_{\tau,h}(t))\| + \| \mvec v(t)- \mvec v_{\tau,h} (t)\| 
	& + \| p(t) - p_{\tau,h}(t)\|\big\}  \leq c (\tau^{k+1}+ h^{r+1})
	\end{aligned}
\end{equation*}
for $t\in I$ is proved. This inequality bounds the error of the total energy \eqref{Eq:HPSenergy}. By Poincar\'{e}'s inequality, an error estimation for $\mvec u$ in $L^\infty(L^2)$ is also induced. However, the resulting $L^\infty(L^2)$ estimate for $\mvec u$ and the bound for $\mvec v$ in \eqref{MainEstEng} are suboptimal only with respect to the discretization in space. For this we recall that $\mvec u$ and $\mvec v$ are approximated in space by piecewise polynomials of order $r+1$. This suboptimality of the error analysis \cite{MB_BAKR22} is linked to the mixed approximation of the Navier--Stokes equations, where the same defect arises. Nevertheless, numerical experiments indicate optimal order convergence in the $L^2(L^2)$ norm; cf.~\cite[Thm.~7.35, Exp.~7.41]{MB_J16}.

In a series of works \cite{MB_ML92,MB_ML94,MB_MTL96}, their authors investigate stability, singular behavior for $t\rightarrow 0$ and decay of solutions to the quasi-static Biot system    
\begin{equation}
\label{Eq:QSBiot}
- \nabla \cdot (\mvec C \mvec \varepsilon (\mvec u)) + \alpha \mvec \nabla p = \mvec f\,, \quad \nabla \cdot \partial_t \mvec u  -  \nabla \cdot (\mvec K \nabla p)  = g
\end{equation}
with the inital condition $\nabla \cdot \mvec u(0) =0$. Further, its approximation by equal order and nonequal order (Taylor--Hood) pairs of finite element spaces is addressed. The observation in \cite{MB_MTL96} is that a behavior of solutions for small $t$ similar to the Navier--Stokes equations is present. They are studied in a sequence of papers by Heywood and Rannacher \cite{MB_HR82}. The conclusion in  \cite{MB_MTL96} is that the consolidation process causes a regularization of exact solutions and  stabilization of the pore pressure field in the unstable (in the sense of inf-sup stability) equal order approximation. Consequently, possible oscillations in the pore pressure field close to $t = 0$ decay in time, particularly for unstable methods. After some time, both stable and unstable methods converge; indeed the equal order approximation is natural and convergent of optimal order.

Driven by \cite{MB_ML92,MB_ML94,MB_MTL96}, here we prove an optimal order error estimate in $L^2(L^2)$ for the approximation of the dynamic model of poroelasticity \eqref{Eq:HPS} by a familiy of continuous space-time finite element methods using an equal order approximation in space of the fields $\mvec u,\mvec v$ and the pore pressure $p$. Thereby, the before-given energy error estimate is sharpened. Computational efficiency is increased by the equal order techniques. The result is illustrated by a numerical experiment.

\section{Numerical scheme}
\label{sec:numerscheme}

For the time discretization, we decompose $I:=(0,T]$ into $N$ subintervals $I_n=(t_{n-1},t_n]$, $n=1,\ldots,N$, where $0=t_0<t_1< \cdots < t_{N-1} < t_N = T$ such that $I=\bigcup_{n=1}^N I_n$. We put $\tau := \max_{n=1,\ldots, N} \tau_n$ with $\tau_n = t_n-t_{n-1}$. The set $\mathcal{M}_\tau := \{I_1,\ldots,
I_N\}$  is called the time mesh. For a Banach space $B$, any $k\in \N$ and $\mathbb P_k(I_n;B) := \big\{w_\tau \,: \,  I_n \to B \,, \; w_\tau(t) = \sum_{j=0}^k  W^j t^j \; \forall t\in I_n\,, \; W^j \in B\; \forall j \big\}$ we define the space of globally continuous, piecewise polynomial functions in time with values in $B$ by 
\begin{equation}
	\label{Eq:DefXk}
	X_\tau^{k} (B) := \big\{w_\tau \in C(\overline I;B) \mid w_\tau{}_{|I_n} \in
	\mathbb P_{k}(I_n;B)\; \forall I_n\in \mathcal{M}_\tau \big\}\,.
\end{equation}
For time integration, it is natural to use the $(k+1)$-point Gau{ss}--Lobatto quadrature formula. On $I_n$, it reads as $Q_n(w) := \frac{\tau_n}{2}\sum_{\mu=1}^{k+1} \hat \omega_\mu w(t_{n,\mu} ) \approx \int_{I_n} w(t) \ud t$,  where $t_{n,\mu} =T_n(\hat t_{\mu})$, for $\mu = 1,\ldots,k+1$, are the quadrature  points on $I_n$ and $\hat \omega_\mu$ the corresponding weights of the quadrature formula. Here, $T_n(\hat t):=(t_{n-1}+t_n)/2 + (\tau_n/2)\hat t$ is the affine transformation from $\hat I = [-1,1]$ to $I_n$ and $\hat t_{\mu}$ are the quadrature points on $\hat I$. The  Gau{ss}--Lobatto formula is exact for all $w\in \mathbb P_{2k-1} (I_n;\R)$. 

For the space discretization, let $\mathcal{T}_h=\{K\}$ be the quasi-uniform decomposition of $\Omega$ into (open) quadrilaterals or hexahedrals, with mesh size $h>0$. Quadrilaterals are used in our numerical experiment of Sec.~\ref{Sec:NumExp} that is based on an inplementation in the deal.II library \cite{MB_Aetal22}. For some integer $r\in \N$, the finite element spaces for approximating the unknowns $\mvec u$, $\mvec v$ and $p$ of \eqref{Eq:HPS} in space are of the form 
\begin{subequations}
\label{Def:VhQh}
\begin{alignat}{2}
\mvec V_{h}^r  & := \{\mvec v_h \in C(\overline \Omega )^d \mid \mvec v_{h}{}_{|K}\in 
{\mvec V^r(K)} \;\; \text{for all}\; K \in \mathcal{T}_h\} \cap \mvec H^1_0(\Omega)\,, \\
	Q_{h}^r  & := \{ q_h \in C(\overline \Omega )\mid q_{h}{}_{|K}\in {Q^r(K)} \;\; \text{for all}\; K \in \mathcal{T}_h\} \cap H^1_0(\Omega)\,.
\end{alignat}
\end{subequations}
For the local spaces $\{\mvec V^r(K),Q^r(K)\}$ we use mapped versions of $\{\mathbb Q_r^d,\mathbb Q_r\}$; cf.\ \cite{MB_J16}.

We use a temporal test basis that is supported on the subintervals $I_n$. Then, a time marching process is obtained. In that, we assume that the trajectories $\mvec u_{ \tau,h}$, $\mvec v_{ \tau,h}$ and $p_{ \tau,h}$ have been computed before for all $t\in [0,t_{n-1}]$, starting with approximations  $\mvec u_{\tau,h}(t_0) :=\mvec u_{0,h}$, $\mvec v_{\tau,h}(t_0) :=\mvec v_{0,h}$ and $p_{\tau,h}(t_0) := p_{0,h}$ of the initial values $\mvec u_0$, $\mvec v_0$ and $p_{0}$. We then consider solving the following local problem on $I_n$; cf.\ \cite{MB_BAKR22}. 

\begin{problem}[$I_n$-problem of continuous, equal order in space approximation]
	\label{Prob:CG}
	Let $k,r\in \N$. For given $\mvec u_{h}^{n-1}:= \mvec u_{\tau,h}(t_{n-1})\in \mvec V_h^r$, $\mvec v_{h}^{n-1}:=  \mvec v_{\tau,h}(t_{n-1})\in \mvec V_h^r$,  and $p_{h}^{n-1}:= p_{\tau,h}(t_{n-1}) \in Q_h^r$ with  $\mvec u_{\tau,h}(t_0) :=\mvec u_{0,h}$, $\mvec v_{\tau,h}(t_0) :=\mvec v_{0,h}$ and $p_{\tau,h}(t_0) := p_{0,h}$, find $(\mvec u_{\tau,h},\mvec v_{\tau,h},p_{\tau,h}) \in \mathbb P_k (I_n;\mvec V_h^r) \times \mathbb P_k (I_n;\mvec V_h^r) \times \mathbb P_k (I_n;Q_h^r)$ such that $\mvec u_{\tau,h}^+(t_{n-1})=\mvec u_{h}^{n-1}$, $\mvec v_{\tau,h}^+(t_{n-1})=\mvec v_{h}^{n-1}$ and $p_{\tau,h}^+(t_{n-1})= p_{h}^{n-1}$ and, for all $(\mvec \phi_{\tau,h},\mvec \chi_{\tau,h},\psi_{\tau,h})\in  \mathbb P_{k-1} (I_n;\mvec V^r_h) \times \mathbb P_{k-1} (I_n;\mvec V^r_h) \times \mathbb P_{k-1} (I_n;Q^r_h)$, 
	\begin{align*}
				Q_n \big(\langle \partial_t \mvec u_{\tau,h} , \mvec \phi_{\tau,h} \rangle  \! - \! \langle \mvec v_{\tau,h} , \mvec \phi_{\tau,h} \rangle \big)  & = 0 \,, \\[1ex]
				Q_n \big(\langle \rho \partial_t \mvec v_{\tau,h} , \mvec \chi_{\tau,h} \rangle \! + \! \langle \mvec C \mvec \varepsilon (\mvec v_{\tau,h}),\mvec \varepsilon (\mvec \chi_{\tau,h})\rangle \! - \! \alpha \langle \nabla \cdot \mvec \chi_{\tau,h}, p_{\tau,h}\rangle\big)  & = Q_n \big(\langle \mvec f , \mvec \chi_{\tau,h} \rangle\big) \,,\\[1ex]
				Q_n \big(\langle c_0 \partial_t p_{\tau,h},\psi_{\tau,h} \rangle \! + \! \alpha \langle \nabla \cdot \mvec v_{\tau,h}, \psi_{\tau,h}\rangle \! + \! \langle \mvec K \nabla p_{\tau,h}, \nabla \psi_{\tau,h} \rangle \big) & = Q_n \big(\langle g,\psi_h \rangle\big)\,.
	\end{align*}
	
\end{problem}

Here, $w^+(t_{n-1}) := \lim_{t \searrow t_{n-1}} w(t)$. The trajectories defined by Problem~\ref{Prob:CG} satisfy that $\mvec u_{\tau,h},\mvec v_{\tau,h}\in X_\tau^k(\mvec V_h^r)$ and $p_{\tau,h}\in X_\tau^k(Q_h^r)$; cf.\ \eqref{Eq:DefXk} and \eqref{Def:VhQh}. Well-posedness of Problem \ref{Prob:CG} can be shown by the arguments of \cite[Lem.\ 3.2]{MB_ABMS23}. Precisely, the scheme of Problem~\ref{Prob:CG} represents a Galerkin--Petrov approach since trial and test spaces differ.  

\section{Error estimation}

For the continuous Galerkin method with equal order approximation in space, that is defined in Problem~\ref{Prob:CG}, there holds the following optimal order error estimate.

\begin{theorem}[Error estimate]
\label{Thm:MainRes}	
For the approximation $(\mvec u_{\tau,h},\mvec v_{\tau,h},p_{\tau,h})$ defined by Problem~\ref{Prob:CG} of the sufficiently regular solution $(\mvec u,\mvec v, p)$ to \eqref{Eq:HPS} and a suitable choice of the discrete initial values $(\mvec u_{0,h},\mvec v_{0,h}, p_{0,h})\in \mvec V_h^r\times \mvec V_h^r \times Q_h^r$, with $\mvec u_{0,h} =\mvec R_h \mvec u_0$ for the elliptic projection $\mvec R_h$ of \eqref{Def:vRh}, there holds for $t\in I$ that
\begin{equation}
		\label{Eq:MR01}
		\| \mvec u(t) - \mvec u_{\tau,h} (t) \|  + \| \mvec v(t) - \mvec v_{\tau,h} (t)\|  + \|p(t) - p_{\tau,h}(t)\| \leq c \tau^{k+1}+ c h^{r+1}\,. 
\end{equation}
\end{theorem}

\begin{mproof}
	The proof follows along the lines of \cite[Thm.\ 4.8]{MB_BAKR22} in that the energy error estimate presented in  Sec.~\ref{Sec:Intro} of this work for inf-sup stable pairs of finite element spaces is shown. However, one of the auxiliary estimates, given in \cite[Lem.\ 4.3]{MB_BAKR22}, needs to be sharpened. This is done below in Lem.~\ref{Lem:Intpol}. For further modification of the estimation in \cite{MB_BAKR22} we refer to Rem.~\ref{Rem:ModProof}.
\end{mproof}

In the following, we need some notation. The operator $\mvec R_h: \mvec H^1_0 \mapsto \mvec V_h^{r}$ defines the elliptic 
projection onto $\mvec V_h^{r}$ such that, for $\mvec w\in \mvec H^1_0$, 
\begin{equation}
	\label{Def:vRh}
	\langle  \mvec C \mvec \varepsilon (\mvec R_h \mvec w), \mvec \varepsilon(\mvec \phi_h) \rangle = \langle  \mvec C \mvec \varepsilon (\mvec w), \mvec \varepsilon(\mvec \phi_h) \rangle
\end{equation}
for all $\mvec \phi_h\in \mvec V_h^{r}$. The operator $R_h: H^1_0 \mapsto V_h^r$ defines the elliptic 
projection onto $V_h^r$ such that, for $w\in H^1_0$,
\begin{equation*}
	\label{Def:Rh}
	\langle  \mvec K \nabla R_h w, \nabla \psi_h \rangle = \langle  \mvec K \nabla w, \nabla \psi_h \rangle   
\end{equation*}
for all $\psi_h\in V_h^r$. For the Gau{ss}--Lobatto quadrature points $t_{n,\mu}$, we define the global Lagrange interpolation operator $I_\tau :C^0(\overline I;L^2)\mapsto X^{k}_\tau(B)$ by means of
\begin{equation*}
	\label{Eq:DefLagIntOp}
	I_\tau w(t_{n,\mu}) = w(t_{n,\mu})\,, \quad \mu=1,\ldots,k+1\,, \; n=1,\ldots, N\,.
\end{equation*}
For given $w\in L^2(I;L^2)$, we define the interpolate $\Pi^{k-1}_\tau w\in L^2(I;L^2)$ such that its restriction $\Pi^{k-1}_\tau w{}_{|I_n}\in\mathbb P_{k-1}(I_n;L^2)$, $n=1,\ldots, N$, is determined  by   
\begin{equation}
	\label{Def:Pi}
	\int_{I_n} \langle \Pi^{k-1}_\tau  w , q \rangle \ud t  =
	\int_{I_n} \langle  w , q \rangle \ud t 
	\qquad\forall\, q\in \mathbb P_{k-1}(I_n;L^2) \,.
\end{equation}
We still introduce some special approximation with optimal order properties; cf.~\cite{MB_KM04}.
\begin{defi}[Special approximation $(\mvec w_1,\mvec w_2)$ of $(\mvec u, \partial_t \mvec u)$]
	\label{Def:W}
	Let $\mvec u\in C^1(\overline I;\mvec H^1_0)$ be given. On $I_n=(t_{n-1},t_n]$ we define 
	\begin{equation}
		\label{Eq:DefW}
		\mvec w_1 := I_\tau \Big( \int_{t_{n-1}}^t \mvec w_2(s)\ud s + \mvec R_h \mvec u(t_{n-1})  \Big)\,,
		\qquad\text{where}\qquad  \mvec w_2 := I_\tau (\mvec R_h\partial_t \mvec u) \,.
	\end{equation}
	Further, we put $\mvec w_1(0) := \mvec R_h \mvec u(0)$. 	
\end{defi}

Let $(\mvec w_1,\mvec w_2)$ be given by Def.~\ref{Def:W}. For the solution of \eqref{Eq:HPS} we put $\mvec U := (\mvec u,\mvec v)$. For the solution of Problem \eqref{Prob:CG} we put $\mvec U_{\tau,h} := (\mvec u_{\tau,h},\mvec v_{\tau,h})$. Then, we split the error by 
\begin{subequations}
\begin{alignat}{2}
	\nonumber
	\mvec U - \mvec U_{\tau,h} & = \begin{pmatrix}
		\mvec u - \mvec u_{\tau,h} \\ \mvec v - \mvec v_{\tau,h}  
	\end{pmatrix}
	= \begin{pmatrix}
		\mvec u - \mvec w_1\\ \mvec v - \mvec w_2
	\end{pmatrix}
	+ \begin{pmatrix}
		\mvec w_1 - \mvec u_{\tau,h} \\ \mvec w_2- \mvec v_{\tau,h}	
	\end{pmatrix}\\
	\label{Eq:ErrU_00}
& =: \begin{pmatrix} \mvec \eta_1\\ \mvec \eta_2 \end{pmatrix} + \begin{pmatrix} \mvec E_{\tau,h}^1\\ \mvec E_{\tau,h}^2 \end{pmatrix} = \mvec \eta + \mvec E_{\tau,h}\,,\\[1ex]
	\label{Eq:ErrP_00}
p - p_{\tau,h} & = p - I_\tau R_h p + I_\tau R_h p - p_{\tau,h} =: \omega +  e_{\tau,h}\,.
\end{alignat}
\end{subequations}

Now, we derive an improved order estimate for $\mvec \eta_1$; cf.~\cite[Lem.~4.3]{MB_BAKR22}.
\begin{lemma}[Improved order estimation of $\partial_t \mvec \eta_1$]
\label{Lem:Intpol}
	{Let $k,r\geq 1$}. For $\partial_t \mvec \eta_1$, defined by \eqref{Eq:ErrU_00}, there holds for $\psi_{\tau,h} \in \mathbb P_{k-1}(I_n;V^{r}_h)$ that 
	\begin{equation}
		\label{Eq:Estdiveta}
		\begin{aligned}		
			\left| \int_{I_n} \langle \nabla \cdot \partial_t \mvec \eta_1, \psi_{\tau,h}\rangle \ud t\right| & \leq c \Big(\tau_n^{k+1}\|\partial_t^{k+2} \mvec u\|_{L^2(I_n;\mvec H^1)} + h^{r+1}  \| \mvec u\|_{H^2(I_n;\mvec H^{r+1})}\Big)\\             
			& \qquad \cdot \big(\| \psi_{\tau,h}\|_{L^2(I_n;L^2)} + \| \nabla \Pi_\tau^{k-1} \psi_{\tau,h}\|_{L^2(I_n;L^2)}\big)\,.
		\end{aligned}
	\end{equation}
\end{lemma}

\begin{mproof}
To prove \eqref{Eq:Estdiveta}, from the first of the definitions in \eqref{Eq:DefW} we conclude that 
	\begin{equation}
		\label{Eq:Estdiveta01}
		\mvec \eta_1 = \mvec u - \mvec w_1 = \mvec u - I_\tau \mvec u + I_\tau \mvec u - I_\tau (\mvec R_h \mvec u) - I_\tau 	\int_{t_{n-1}}^t (\mvec w_2 - \partial_t \mvec R_h \mvec u) \ud s\,.
	\end{equation}
By \eqref{Eq:Estdiveta01} we then get that 
\begin{align}
	        \nonumber
			& \int_{I_n} \langle \nabla \cdot \partial_t \mvec \eta_1, \psi_{\tau,h}\rangle \ud t =  \int_{I_n} \langle \nabla \cdot \partial_t (\mvec u - I_\tau \mvec u ), \psi_{\tau,h}\rangle \ud t\\
			\nonumber
			&  \qquad + \int_{I_n} \langle \nabla \cdot \partial_t I_\tau (\mvec u  - \mvec R_h \mvec u ), \psi_{\tau,h}\rangle \ud t\\		
			\label{Eq:Estdiveta02}
			& \qquad  + \int_{I_n} \Big\langle \nabla \cdot \partial_t I_\tau 	\int_{t_{n-1}}^t (\mvec w_2 - \partial_t \mvec R_h \mvec u) \ud s, \psi_{\tau,h}\Big\rangle \ud t   =: \Gamma_1 + \Gamma_2 + \Gamma_3\,.
\end{align}
	
	We start with estimating $\Gamma_1$. Firstly, let $k\geq 2$. Using integration by parts in time and recalling that the endpoints of $I_n$ are included in the set of Gauss--Lobatto quadrature points  of $I_n$, we get that 
	\begin{equation*}
		\Gamma_1 = \int_{I_n} \langle \nabla \cdot \partial_t (\mvec u - I_\tau \mvec u ), \psi_{\tau,h}\rangle \ud t = - \int_{I_n} \langle \nabla \cdot (\mvec u - I_\tau \mvec u), \partial_t \psi_{\tau,h}\rangle \ud t\,.
	\end{equation*}
	Let now $I_\tau^{k+1}$ denote the Lagrange interpolation operator at the $k+2$ points of $\overline I_n = [t_{n-1},t_n]$ consisting of the $k+1$ Gauss--Lobatto quadrature nodes $t_{n,\mu}$, for $\mu=1,\ldots,k+1$, and a further node in $(t_{n-1},t_n)$ that is distinct from the previous ones. Then, $(I_\tau^{k+1}\mvec u)\partial_t \psi_{\tau,h}$ is a polynomial of degree $2k-1$ in $t$, such that 
	\begin{equation*}
		\int_{I_n} \langle \nabla \cdot (\mvec u - I_\tau \mvec u), \partial_t \psi_{\tau,h}\rangle \ud t = \int_{I_n} \langle \nabla \cdot (\mvec u - I_\tau^{k+1} \mvec u), \partial_t \psi_{\tau,h}\rangle \ud t\,.
	\end{equation*}
	Using integration by parts in time and the approximation properties of the Lagrange interpolator $I_\tau^{k+1}$ in the norm of $L^2(I_n;H^1)$, we have that
	\begin{align}
		\nonumber
			|\Gamma_1| & \leq \Big|\int_{I_n} \langle \nabla \cdot \partial_t (\mvec u - I_\tau^{k+1} \mvec u), \psi_{\tau,h}\rangle \ud t\Big|\\[1ex] & 
		    \label{Eq:Estdiveta03}
			\leq  c \tau_n^{k+1}\| \partial_t^{k+2} \mvec u\|_{L^2(I_n;\mvec H^1)} \| \psi_{\tau,h}\|_{L^2(I_n;L^2)}\,.
	\end{align}
	For $k=1$, $(\partial_t I_\tau \mvec u) \psi_{\tau,h}\in \mathbb P_0(I_n;V_h^r)$ and $\partial_t I_\tau \mvec u = (\mvec u(t_n)-\mvec u(t_{n-1}))/\tau_n$ imply that
	\begin{equation}
		\label{Eq:Estdiveta03k1}
		\begin{aligned}
			\Gamma_1 = \Big\langle \nabla \cdot \int_{I_n}  (\partial_t \mvec u - \partial_t I_\tau\mvec u)\ud t , \psi_{\tau,h}\Big\rangle
			=  0\,.
		\end{aligned}
	\end{equation}
	
Next, we estimate $\Gamma_2$. For this we introduce the abbreviation $\mvec \xi := \mvec u -\mvec R_h \mvec u$. The Lagrange interpolant $I_\tau$ satisfies the stability results (cf.~\cite[Eqs.~(3.15) and (3.16)]{MB_KM04})
	\begin{subequations}
		\label{Eq:Estdiveta035}
		\begin{alignat}{2}
			\label{Eq:Estdiveta0351}
			\|  I_\tau w \|_{L^2(I_n;L^2)} & \leq c \| w \|_{L^2(I_n;L^2)}  + c \tau_n \|  \partial_t w \|_{L^2(I_n;L^2)} \,, \\[1ex]
			\label{Eq:Estdiveta0352}
			\Big\|  \int_{t_{n-1}}^t w \ud s \Big \|_{L^2(I_n;L^2)}  & \leq c \tau_n  \|  w \|_{L^2(I_n;L^2)} \,.
		\end{alignat}
	\end{subequations}
Viewing $\mvec \xi(t_{n-1}^+)$ as a function constant in time, using integration by parts for the space variables and recalling the definition \eqref{Def:Pi} of $\Pi_\tau^{k-1}$, we find that 
\begin{equation*}
	\label{Eq:Estdiveta04}
	\begin{aligned}
		\Gamma_2  &  = \int_{I_n} \langle \nabla \cdot \partial_t I_\tau (\mvec \xi- \mvec \xi(t_{n-1}^+)), \psi_{\tau,h}\rangle \ud t = \int_{I_n} \langle \nabla \cdot \partial_t I_\tau \int_{t_{n-1}}^t \partial_t \mvec \xi\ud s, \psi_{\tau,h}\rangle \ud t\\[1ex] 
		& = \int_{I_n} \langle \partial_t I_\tau \int_{t_{n-1}}^t \partial_t \mvec \xi\ud s, \nabla \psi_{\tau,h}\rangle \ud t= \int_{I_n} \langle \partial_t I_\tau \int_{t_{n-1}}^t \partial_t \mvec \xi\ud s, \nabla \Pi_\tau^{k-1}\psi_{\tau,h}\rangle \ud t\,.
	\end{aligned}
\end{equation*}
By the $H^1$--$L^2$ inverse inequality $\| w'\|_{L^2(I_n;\R)}\leq c \tau_n^{-1}\| w\|_{L^2(I_n;\R)}$, the stability results \eqref{Eq:Estdiveta04}, the standard approximation properties of the elliptic projection $\mvec R_h$ defined in \eqref{Def:vRh}, the latter identity yields that
\begin{align}
	\nonumber
			|\Gamma_2 | &  = \Big| \int_{I_n} \langle \partial_t I_\tau \int_{t_{n-1}}^t \partial_t \mvec \xi\ud s, \nabla \Pi_\tau^{k-1}\psi_{\tau,h}\rangle \ud t\Big|\\[1ex] 
	\nonumber
			& \leq c \tau_n^{-1} \Big\| I_\tau \int_{t_{n-1}}^{t} \partial_t \mvec \xi \ud s
			\Big\|_{L^2(I_n;L^2)} \| \nabla \Pi_\tau^{k-1}\psi_{\tau,h}\|_{L^2(I_n;L^2)}\\[1ex]
		    \label{Eq:Estdiveta04}
			& \leq c h^{r+1}\| \partial_t \mvec u \|_{L^2(I_n;\mvec H^{r+1})}\| \nabla \Pi_\tau^{k-1}\psi_{\tau,h}\|_{L^2(I_n;L^2)}\,.
\end{align}

Finally, we estimate $\Gamma_3$. By the arguments of \eqref{Eq:Estdiveta04} it follows for $\Gamma_3$ that 
\begin{align*}
			|\Gamma_3 | & = \Big| \int_{I_n} \Big\langle \partial_t I_\tau \int_{t_{n-1}}^t (\mvec w_2 - \partial_t \mvec R_h \mvec u) \ud s, \nabla \Pi_\tau^{k-1} \psi_{\tau,h}\Big\rangle \ud t\Big| \\[1ex]
			& \leq c \|\mvec w_2  - \mvec R_h (\partial_t \mvec u)  \|_{L^2(I_n;\mvec L^2)}\|\nabla \Pi_\tau^{k-1} \psi_{\tau,h}\|_{L^2(I_n;L^2)}\,.
\end{align*}
By \cite[Lem.~3.3]{MB_KM04}, we then obtain that
	\begin{equation}
		\label{Eq:Estdiveta05}
		\begin{aligned}
			|\Gamma_3 | & \leq c \big(\tau_n^{k+1}\|\partial_t^{k+2} \mvec u \|_{L^2(I_n;\mvec L^2)} + h^{r+1}\| \mvec u \|_{H^2(I_n;\mvec H^{r+1})}\big)\\
			& \quad  \cdot  \|\nabla \Pi_\tau^{k-1} \psi_{\tau,h}\|_{L^2(I_n;L^2)}\,.
		\end{aligned}
	\end{equation} 
Combining \eqref{Eq:Estdiveta02} with \eqref{Eq:Estdiveta03}, \eqref{Eq:Estdiveta03k1}, \eqref{Eq:Estdiveta04} and \eqref{Eq:Estdiveta05} proves the assertion \eqref{Eq:Estdiveta}. 
\end{mproof}

For the proof of \eqref{Eq:MR01} we note the following. 
\begin{rem}
	\label{Rem:ModProof}
	\begin{itemize}
		\item In the stability estimate for $\mvec E_{\tau,h}$ and $e_{\tau,h}$ of the equal order scheme of Problem~\ref{Prob:CG}, replacing the stability result of \cite[Lem.~4.4]{MB_BAKR22} for the Taylor--Hood finite element pair based approach, the term  $\|\nabla \Pi_\tau^{k-1} \psi_{\tau,h}\|_{L^2(I_n;L^2)}$ arising now in \eqref{Eq:Estdiveta} is absorbed by $\int_{I_n} \langle K \nabla \Pi_\tau^{k-1}  e_{\tau,h}, \nabla \Pi_\tau^{k-1}  e_{\tau,h},  \rangle  \ud t$ on the left-hand side of the stability inequality for $\mvec E_{\tau,h}$ and $e_{\tau,h}$. For this, the inequalities of Cauchy--Schwarz and Cauchy--Young have to be applied to the right-hand side of \eqref{Eq:Estdiveta}.
		
		\item The choice $\mvec u_{0,h} =\mvec R_h \mvec u_0$ of $\mvec u_{0,h}$ is essential for the analysis. In computations, the sharpness of this assumption is not observed; cf.~\cite{MB_BAKR22}, \cite[Rem.~3.1]{MB_KM04}. The discrete initial values have to ensure that the convergence rate is not perturbed. 
		
		\item In \cite{MB_BAKR22}, the proof of the error estimate is not build on an inf-sup stability condition for the involved spatial finite element spaces which allows its adaptation to the equal order approximation studied here. Limit cases of vanishing coefficients in \eqref{Eq:HPS}, for instance incompressible material, are not covered by Thm.~\ref{Thm:MainRes}. They do not become feasible without further assumptions like the inf-sup stability condition. 
	\end{itemize}
\end{rem}

\section{Numerical convergence study}
\label{Sec:NumExp}

Here we present the results of our performed numerical experiment to illustrate Thm.~\ref{Thm:MainRes}.
The implementation was done in an in-house frontend solver for the \texttt{deal.II} library \cite{MB_Aetal22}. We study \eqref{Eq:HPS} for $\Omega=(0,1)^2$ and $I=(0,2]$ and the prescribed solution
\begin{equation}
	\label{Eq:givensolution}
	\begin{aligned}
	& \mvec u(\boldsymbol x, t) = \phi(\mvec x, t) \boldsymbol I_2 \;\; \text{and}\;\;
	p(\mvec x, t) = \phi(\mvec x, t)\\
	& \text{with}\;\;  \phi(\mvec x, t) = \sin(\omega_1 t^2) \sin(\omega_2 x_1) \sin(\omega_2 x_2)
	\end{aligned}
\end{equation}
and $\omega_1=\omega_2 = \pi$. We put $\rho=1.0$, $\alpha=0.9$, $c_0=10^{-3}$ and $\boldsymbol K= 10^{-2} \cdot \boldsymbol I_2$ with the identity $\mvec I_2\in \R^{2,2}$. For the fourth order elasticity tensor $\boldsymbol C$, isotropic material properties with Young's modulus $E=100$ and Poisson's ratio $\nu=0.35$ are chosen. The norm of $L^\infty(I;L^2)$ is approximated by $	\| w\|_{L^\infty(I;L^2)} \approx \max \{ \| w_{|I_n}(t_{n,m})\|  \mid m=1,\ldots ,M\,, \; n=1,\ldots,N\}$ with $M=100$ and the Gauss quadrature nodes $t_{n,m}$ of $I_n$. We investigate the space-time convergence of the scheme in Problem~\ref{Prob:CG} to confirm \eqref{Eq:MR01}. For this, the domain $\Omega$ is decomposed into a sequence of successively refined meshes of quadrilateral finite elements. The spatial and temporal mesh sizes are halfened in each of the refinement steps. The step sizes of the coarsest mesh are $h_0=1/(2\sqrt 2)$ and $\tau_0=0.1$. We choose the polynomial degree $k=2$ and $r=2$, such that a solution $(\mvec u_{\tau,h}, \mvec v_{\tau,h})\in (X_\tau^2(\mvec V_h^2))^2$, $p_{\tau,h}\in X_\tau^2(Q_h^2)$ is obtained; cf.\ \eqref{Eq:DefXk} and  \eqref{Def:VhQh}. The calculated errors and corresponding experimental orders of convergence are summarized in Tab.~\ref{Tab:1}. Tab.~\ref{Tab:1} confirms our main result \eqref{Eq:MR01}. A severe lack of stability of the equal order in space approximation is not observed. Optimal order approximation properties are ensured. For the pressure variable the asymptotic range of convergence is not completely reached yet (compared to \cite[Sec.~5]{MB_BAKR22}), which is supposed to be due to the reduced stability of the equal order in space approach. 

\begin{table}[H]
	\centering
\begin{tabular}{l}
	\begin{tabular}{cccccccc}
		\toprule
		{$\tau$} & {$h$} &
		{ $\| \mvec u - \mvec u_{\tau,h}  \|_{L^2(\mvec L^2)} $ } & {EOC} &
		{ $\| \mvec v - \mvec v_{\tau,h}  \|_{L^2(\mvec L^2)} $ } & {EOC} &
		{ $\| p-p_{\tau,h}  \|_{L^2(L^2)}  $ } & {EOC} 
		\\
		\cmidrule(r){1-2}
		\cmidrule(lr){3-8}
		$\tau_0/2^0$ & $h_0/2^0$ & 3.6741729210e-03 & {--} & 4.1282666037e-02 & {--} & 6.7442978958e-02 & {--}  \\ 
		$\tau_0/2^1$ & $h_0/2^1$ & 4.6822546563e-04 & 2.97 & 4.9051780039e-03 & 3.07 & 3.8674819986e-03 & 4.12  \\
		$\tau_0/2^2$ & $h_0/2^2$ & 5.8778097800e-05 & 2.99 & 6.2028252749e-04 & 2.98 & 2.4809531869e-04 & 3.96  \\
		$\tau_0/2^3$ & $h_0/2^3$ & 7.3543857134e-06 & 3.00 & 7.7854545813e-05 & 2.99 & 1.6114838644e-05 & 3.94  \\
		$\tau_0/2^4$ & $h_0/2^4$ & 9.1951071480e-07 & 3.00 & 9.7424834053e-06 & 3.00 & 1.1541868101e-06 & 3.80  \\
		$\tau_0/2^5$ & $h_0/2^5$ & 1.1494533804e-07 & 3.00 & 1.2181487408e-06 & 3.00 & 1.0079612177e-07 & 3.52  \\
		\bottomrule
	\end{tabular}\\
	\mbox{}\\
	\begin{tabular}{cccccccc}
		\toprule
		{$\tau$} & {$h$} &
		{ $\| \mvec u - \mvec u_{\tau,h}   \|_{L^{\infty}(\mvec L^2)} $ } & {EOC} &
		{ $\| \mvec v - \mvec v_{\tau,h}  \|_{L^{\infty}(\mvec L^2)} $ } & {EOC} &
		{ $\| p-p_{\tau,h}  \|_{L^{\infty}(L^2)} $ } & {EOC}
		\\
		\cmidrule(r){1-2}
		\cmidrule(lr){3-8}
		$\tau_0/2^0$ & $h_0/2^0$ & 7.7356665842e-03 & {--} & 1.1023905122e-01 & {--} & 1.4136586176e-01 & {--}  \\ 
		$\tau_0/2^1$ & $h_0/2^1$ & 1.2393232164e-03 & 2.64 & 1.4136203260e-02 & 2.96 & 8.3752373581e-03 & 4.08  \\
		$\tau_0/2^2$ & $h_0/2^2$ & 1.6642245699e-04 & 2.90 & 1.8234611143e-03 & 2.95 & 4.9215141823e-04 & 4.09  \\
		$\tau_0/2^3$ & $h_0/2^3$ & 2.1296270280e-05 & 2.97 & 2.2979190855e-04 & 2.99 & 2.8859565155e-05 & 4.09  \\
		$\tau_0/2^4$ & $h_0/2^4$ & 2.6851742263e-06 & 2.99 & 2.8663447685e-05 & 3.00 & 2.4527886661e-06 & 3.56  \\
		$\tau_0/2^5$ & $h_0/2^5$ & 3.3684383499e-07 & 2.99 & 3.5800946159e-06 & 3.00 & 2.5725548188e-07 & 3.25  \\
		\bottomrule
	\end{tabular}
\end{tabular}
	\caption{%
		$L^2(L^2)$ and $L^\infty(L^2)$ errors and experimental orders of convergence (EOC) for \eqref{Eq:givensolution} and the equal order scheme of Problem~\ref{Prob:CG} with polynomial degree  $k=2$ and $r=2$ in \eqref{Eq:DefXk} and  \eqref{Def:VhQh}.
	}
	\label{Tab:1}
\end{table}

For the sake of completeness, we still summarize in Tab.~\ref{Tab:2} the computational results that we obtained for the inf-sup stable pair $\{\mvec V_h^3,Q_h^2\}$ of Taylor--Hood finite element spaces (cf.~\cite{MB_Z22}) leading to an approximation $(\mvec u_{\tau,h}, \mvec v_{\tau,h})\in (X_\tau^2(\mvec V_h^3))^2$ and $p_{\tau,h}\in X_\tau^2(Q_h^2)$ of the solution to \eqref{Eq:HPS}; cf.~\cite{MB_ABMS23} for the application of inf-sup stable schemes to \eqref{Eq:HPS}. No striking difference to the computationally less expensive equal order approximation of Tab.~\ref{Tab:1} is observed. This confirms our conjecture of the introduction that the equal-order approximation in space of  \eqref{Eq:HPS} is natural and convergent of optimal order. For the incompressible limit with $c_0 \rightarrow 0$, $\mvec K\rightarrow \mvec 0$ and $g\equiv 0$, robustness of the convergence cannot be expected for the equal-order scheme, since a system of Stokes-type structure with the well-known stability issues of mixed approximations is obtained.  

\begin{table}[H]
	\centering
	\begin{tabular}{l}
		\begin{tabular}{cccccccc}
			\toprule
			{$\tau$} & {$h$} &
			{ $\| \mvec u - \mvec u_{\tau,h}  \|_{L^2(\mvec L^2)} $ } & {EOC} &
			{ $\| \mvec v - \mvec v_{\tau,h}  \|_{L^2(\mvec L^2)} $ } & {EOC} &
			{ $\| p-p_{\tau,h}  \|_{L^2(L^2)}  $ } & {EOC} 
			\\
			\cmidrule(r){1-2}
			\cmidrule(lr){3-8}
			$\tau_0/2^0$ & $h_0/2^0$ & 3.2223164159e-03 & {--} & 3.9470120875e-02 & {--} & 5.3817187010e-02 & {--}  \\ 
			$\tau_0/2^1$ & $h_0/2^1$ & 3.8762468516e-04 & 3.06 & 4.4805250509e-03 & 3.14 & 3.5640199518e-03 & 3.92  \\
			$\tau_0/2^2$ & $h_0/2^2$ & 4.8084739364e-05 & 3.01 & 5.6306412015e-04 & 2.99 & 3.0169432186e-04 & 3.56  \\
			$\tau_0/2^3$ & $h_0/2^3$ & 6.0009965689e-06 & 3.00 & 7.0594365429e-05 & 3.00 & 3.1449270169e-05 & 3.26  \\
			$\tau_0/2^4$ & $h_0/2^4$ & 7.4984738840e-07 & 3.00 & 8.8318214147e-06 & 3.00 & 3.7139677858e-06 & 3.08  \\
			\bottomrule
		\end{tabular}\\
		\mbox{}\\
		\begin{tabular}{cccccccc}
			\toprule
			{$\tau$} & {$h$} &
			{ $\| \mvec u - \mvec u_{\tau,h}   \|_{L^{\infty}(\mvec L^2)} $ } & {EOC} &
			{ $\| \mvec v - \mvec v_{\tau,h}  \|_{L^{\infty}(\mvec L^2)} $ } & {EOC} &
			{ $\| p-p_{\tau,h}  \|_{L^{\infty}(L^2)} $ } & {EOC}
			\\
			\cmidrule(r){1-2}
			\cmidrule(lr){3-8}
			$\tau_0/2^0$ & $h_0/2^0$ & 7.9770740804e-03 & {--} & 1.1179102428e-01 & {--} & 1.6306452013e-01 & {--}  \\ 
			$\tau_0/2^1$ & $h_0/2^1$ & 1.2538478867e-03 & 2.67 & 1.4112728641e-02 & 2.99 & 1.2187980808e-02 & 3.74  \\
			$\tau_0/2^2$ & $h_0/2^2$ & 1.6721629610e-04 & 2.91 & 1.8029843969e-03 & 2.97 & 7.2946116772e-04 & 4.06  \\
			$\tau_0/2^3$ & $h_0/2^3$ & 2.1341460866e-05 & 2.97 & 2.2791756297e-04 & 2.98 & 3.8469033488e-05 & 4.25  \\
			$\tau_0/2^4$ & $h_0/2^4$ & 2.6878519380e-06 & 2.99 & 2.8496047406e-05 & 3.00 & 4.0634844335e-06 & 3.24  \\
			\bottomrule
		\end{tabular}
	\end{tabular}
	\caption{%
		$L^2(L^2)$ and $L^\infty(L^2)$ errors and experimental orders of convergence (EOC) for solution \eqref{Eq:givensolution}, the inf-sup stable pair $\{\mvec V_h^3,Q_h^2\}$ of Taylor--Hood finite element spaces for the approximation in space and polynomial degree  $k=2$  in \eqref{Eq:DefXk}  for the  approximation in time. 
	}
	\label{Tab:2}
\end{table}

\section*{Acknowledgment}
	Computational resources (HPC-cluster HSUper) have been provided by the project hpc.bw, funded by dtec.bw --- Digitalization and Technology Research Center of the Bundeswehr. dtec.bw is funded by the European Union --- NextGenerationEU.


\begin{thebibliography}{99.}%


\bibitem{MB_Aetal22}
Arndt, D.\  et al.: {The deal.II Library, Version 9.4}. J.\ Numer.\ Math.\ \textbf{30}, 231--245 (2022)

\bibitem{MB_ABMS23}
Anselmann, M., Bause, M., Margenberg, N., Shamko, P.: {An energy-efficient GMRES--Multigrid solver for space-time finite element computation of dynamic poro- and thermoelasticity}. Comput.\ Mech.\ \textbf{submitted}, 1--30 (2023); arXiv:2303.06742
 
\bibitem{MB_BAKR22}
Bause, M., Anselmann, M., Köcher, U., Radu, F.A.: {Convergence of a continuous Galerkin method for hyperbolic-parabolic systems}.  Comput.\ Math.\ with Appl.\ \textbf{158}, 118--138 (2024); doi: 10.1016/j.camwa.2024.01.014

\bibitem{MB_BF91}
Brezzi, F., Falk, R.: {Stability of higher-order Hood-Taylor methods}. SIAM J.\ Numer.\ Anal.\ \textbf{28}, 581--590 (1991)

\bibitem{MB_HR82} 
Heywood, J., Rannacher, R.: {Finite element approximation of the nonstationary Navier--Stokes problem. I. Regularity of solutions and second-order error estimates for spatial discretization}. SIAM J.~Numer.~Anal.\ \textbf{19}, 275--311 (1982)
	
\bibitem{MB_JR18}
Jiang, S., Racke, R.: {Evolution equations in thermoelasticity}. CRC Press, Boca Raton (2018)

\bibitem{MB_J16}
John, V.: {Finite Element Methods for Incompressible Flow Problems}. Springer, Cham (2016)

\bibitem{MB_KM04}
Karakashian, O., Makridakis, C.: {Convergence of a continuous Galerkin method with mesh modification for nonlinear wave equations}. Math.\ Comp.\ \textbf{74}, 85--102 (2004)

\bibitem{MB_ML92}
Murad, M.A., Loula, A.F.D.: {Improved accuracy in finite element analysis of Biot’s consolidation problem}. Comput.\ Methods Appl.\ Mech. Engrg.\ \textbf{95}, 359--382 (1992)

\bibitem{MB_ML94}
Murad, M.A., Loula, A.F.D.: {On stability and convergence of finite element approximations of Biot’s consolidation problem}. Int.\ J.\ Numer.\ Methods Eng.\ \textbf{37}, 645--667 (1994)

\bibitem{MB_MTL96}
Murad, M.A., Thom\'{e}e, V., Loula, A.F.D.: {Asymptotic behavior of semidiscrete finite-element approximations of Biot’s consolidation problem}. SIAM J.\ Numer.\ Anal.\ \textbf{33}, 1065--1083 (1996)
	
\bibitem{MB_STW22}
Seifert, C., Trostorff, S., Waurick, M.: {Evolutionary Equations: Picard's Theorem for Partial Differential Equations, and Applications}. Birkhäuser, Cham (2022)

\bibitem{MB_S89}
Slodi\v{c}ka, M.: {Application of Rothe's method to integrodifferential equation}. Comment.\ Math.\ Univ.\ Carolinae\ \textbf{30} (1989), 57--70 (1989)

\bibitem{MB_Z22}
W.\ Zulehner, {A short note on inf-sup conditions for the Taylor--Hood family $Q_k$--$Q_{k-1}$}. arXiv, 1--15 (2022); doi: 10.48550/arXiv.2205.14223


\end{thebibliography}
\end{document}